\documentstyle[amsfonts,amssymb,graphicx]{article}
\newtheorem{theorem}{Theorem}[section]
\newtheorem{corollary}[theorem]{Corollary}

\newtheorem{proposition}[theorem]{Proposition}
\newtheorem{conjecture}[theorem]{Conjecture}
\newtheorem{definition}[theorem]{Definition}
\newtheorem{remark}[theorem]{Remark}

\newcommand{\mb}{\mathbb }

\newcommand{\be}{\begin{equation}}
\newcommand{\ee}{\end{equation}}
\newcommand{\beq}{\begin{eqnarray}}
\newcommand{\eeq}{\end{eqnarray}}

\title{The maximum principle and sign changing solutions of the hyperbolic equation with the Higgs potential}

\author{Andras Balogh and Karen Yagdjian}

\date{}

\begin{document}
\maketitle
\thispagestyle{empty}
\maketitle
\thispagestyle{empty}

 \centerline{School of Mathematical and Statistical Sciences,
University of Texas RGV,}
   \centerline{1201 W.~University Drive,
Edinburg, TX 78539,
USA }

\bigskip

\begin{abstract}
In this article we  discuss  the maximum principle for the linear equation and the sign changing solutions of the semilinear equation with the Higgs potential.   
Numerical simulations indicate that the bubbles for  the semilinear Klein-Gordon equation in the  de~Sitter space-time are created and apparently exist  for all times.
\end{abstract}

\section{Introduction}
\label{S1}

 In this article we  discuss  the maximum principle for the linear equation and the sign changing solutions of the semilinear equation with the Higgs potential. 
The Klein-Gordon equation with the Higgs potential  (the  Higgs boson equation) in the
  de~Sitter space-time is the equation
\begin{equation}
\label{HBE}
\psi _{tt} - e^{-2t} \Delta \psi  +n\psi_t  =  \mu^2  \psi  - \lambda \psi  ^3 ,
\end{equation}
where $\Delta  $ is the Laplace operator in $x \in {\mathbb R}^n$, $n=3$, $t >0$,  $\lambda >0 $, and $\mu >0 $. We assume that 
$\psi =\psi (x,t) $ is a real-valued function.

\medskip

We focus on the zeros of the solutions to the linear and semilinear hyperbolic equation in the Minkowski and de~Sitter space-times.
One motivation for the study of the maximum principle,  sign changing solutions and zeros of the solutions to the linear and semilinear hyperbolic equation comes from the cosmological contents and quantum field theory. It is of  considerable interest for  particle physics and inflationary  cosmology  to study the so-called bubbles \cite{Coleman},
\cite{Linde}, \cite{Voronov}.
In \cite{Lee-Wick} bubble is defined as a simply connected domain surrounded by a wall such that the field
 approaches one of the vacuums outside of a bubble. 
 The creation and
growth of bubbles is an interesting mathematical problem \cite[Ch.7]{Coleman},
\cite{Linde}.   In this paper, for the continuous solution $\psi =\psi (x,t)  $ to the Klein-Gordon equation, for  every  given positive time $t$ we  define a bubble  as a maximal connected  set of  points $ x \in {\mathbb R}^n$ at which  solution  changes  sign. 
\smallskip

Another motivation to study all these closely related properties comes from the issue of the existence of a global in time solution  to non-linear equation. Consider the Cauchy problem for the linear wave equation
\begin{eqnarray*}
\cases{  \partial_{t }^2 u- \Delta u=0 \,, \cr
u(x,0)= u_0(x)\,,\quad u_t(x,0)= u_1(x)\,,\qquad u_0,u_1 \in C^\infty ({\mathbb R}^n)\,.}
\end{eqnarray*}
If one can prove that the solution $u=u(x,t) $ vanishes at some point $(x_b,t_b)$, then  it opens the door to study  the  blowup   phenomena for the equation
\begin{eqnarray*}
\cases{  \partial_{t }^2 v- \Delta v + (\partial_{t } v)^2- |\nabla v|^2=0 \,, \cr
v(x,0)= v_0(x)\,,\quad v_t(x,0)= v_1(x)\,,\qquad v_0,v_1 \in C ^\infty ({\mathbb R}^n)\,,}
\end{eqnarray*}
which is  the Nirenberg's Example (see, e.g.\cite{Klainerman1980}) of the quasilinear equation. Indeed, the  
transformation
\[
u(x,t)=\exp(v(t,x))
\]
shows that $u(x_b,t_b) =0$ at some $t_b>0$ and $x_b \in {\mathbb R}^n$ implies $v(x_b,t_b)=-\infty$. To avoid blowup phenomena one can restrict the initial data to be small in some norm.
(For details, see, e.g., \cite{Yag2005}.)
Therefore to guarantee existence of the global solution to quasilinear equation,  the  solution $u=u(x,t)$ of the related linear equation must keep sign for all $t>0$ and all $x  \in {\mathbb R}^n$. This link between sign preserving solutions and global in time solvability is especially easy to trace in the case of $n=3$. In fact, the explicit representation formulas for the solutions to the linear equation play key role. On the other hand for the equations with the variable coefficients and, in particular,  for the linear hyperbolic equations in the curved space-time, the new global in time explicit representation formulas were obtained very recently (see, \cite{Yag_Galst_CMP,MN}).
For the results on   the sign changing solutions of the quasilinear equations one can consult \cite{Speck}.
\medskip

The outline  of the discussion in this   paper is organized as follows. In Section~\ref{S2}  we describe the maximum principle for the wave equation in the Minkowski space-time when the initial data are subharmonic or superharmonic. In   Section~\ref{S3} we present the maximum principle for the 
linear Klein-Gordon equation  in the de~Sitter space-time.  Theorem~\ref{T2} of that section guarantees  that   the solution does not changes
 sign, that is, it provides with  some necessary conditions to have a sign-changing solution. Section~\ref{S5} is devoted to kernels of the integral transforms have been used in the proofs. Section~\ref{S6} is   a bridge between Section~\ref{S7} and previous sections. It    is aimed to give some theoretical background material about semilinear Klein-Gordon equation  in the de~Sitter space-time with the Higgs potential. 
Section~\ref{S6} also prepares the reader to  Section~\ref{S7}, which is about numerical simulations on the evolution of the bubbles in the de~Sitter space-time.

\section{The maximum principle in the Minkowski space-time}
\label{S2}

In \cite{Sather} the following maximum principle is established
for   the wave operator
\[
L:=\partial_{t }^2- \Delta\,,
\]
where $\Delta  $ is the Laplace operator in  $x \in {\mathbb R}^n $. Denote
\begin{eqnarray*}
N= \cases{ \frac{n-2}{2}\quad if \quad n \quad even \cr
\frac{n-3}{2}\quad if \quad n \quad  odd.}
\end{eqnarray*}
Let   $ u$ satisfy the differential inequality
\begin{eqnarray*}
\frac{ \partial ^N}{\partial t ^N}  (L [ u] ) \leq 0    \quad \mbox{\rm for all} \quad t \in [0, T]\,,
\end{eqnarray*}
and the initial conditions
\begin{eqnarray*}
&  &
\frac{ \partial ^k u}{\partial t ^k} (x,0)=0\,, \quad  k=0,1,\ldots,N\,,   \quad
\frac{ \partial ^{N+1}u}{\partial t ^{N+1}} (x,0)\leq 0\,,
\end{eqnarray*}
for all $x$ in the  domain $D_0\subseteq  {\mathbb R}^n $. Then
\[
u (x,t) \leq 0\,
\]
in the domain of dependence of $D_0$, where $t\leq T $.

We recall definition of the  forward light cone  $D_+ (x_0,t_0) $,
and   the  backward light cone $D_- (x_0,t_0) $, in the Minkowski space-time  for the point  $(x_0, t_0) \in {\mathbb R}^{n+1} $:
\begin{eqnarray*}
D_\pm (x_0,t_0)
& := &
\Big\{ (x,t)  \in {\mathbb R}^{n+1}  \, ; \,
|x -x_0 | \leq \pm(  t-t_0 )
\,\Big\} \,.
\end{eqnarray*}
For the domain $D_0\subseteq {\mathbb R}^n$ define a dependence domain of $D_0 $ as follows:
\begin{eqnarray*}
D(D_0)
& := &
\displaystyle  \bigcup_{x_0 \in {\mathbb R}^n,\,\, t_0 \in [0,\infty) }    \left\{ D_- (x_0,t_0) \,;\,  D_- (x_0,t_0) \cap \{t=0\} \subset D_0 \right\}\,.
\end{eqnarray*}
In particular, according to Theorem~1~\cite{Sather},   for   $x \in {\mathbb R}^3 $ if $ u$ satisfy the differential inequality
\begin{eqnarray}
\label{2}
L  u  \leq 0 \,,  \quad \mbox{\rm for all} \quad t \leq T\,,
\end{eqnarray}
and the initial conditions
\[
u(x,0)=0\,, \quad  u_t(x,0) \leq 0\,,  \quad \mbox{\rm for all} \quad x \in D_0 \subseteq {\mathbb R}^3\,,
\]
then $ u (x,t) \leq 0$ in the domain of dependence of $D_0$, where $t\leq T $. The statement is a simple consequence of the Duhamel's  principle and the well-known Kirchhoff  formula
\begin{eqnarray}
\label{6}
u(x_1,x_2,x_3,t)= \frac{1}{4\pi}\int_{S_t(x_1,x_2,x_3)} \frac{\varphi (\alpha _1,\alpha _2,\alpha _3)}{t} d S_t \,,
\end{eqnarray}
for the solution of the Cauchy problem for the wave equation (see, e.g., \cite{Shatah}), where
\[
L[u]=0\,, \qquad u(x,0)=  0\,,\qquad u_t(x,0)= \varphi (x)\,,
\]
and $ S_t(x_1,x_2,x_3)$ is a sphere of radius $t$ centered at $(x_1,x_2,x_3) $. The kernel $1/t$ of the integral operator   (\ref{6}) is positive.
The next statement  also can be proved by the Kirchhoff formula. 

\begin{theorem}
\label{T1.1b}
Assume that the  $C^2$- functions $u=u(x,t)$, $ \varphi _0= \varphi _0(x)$, $\varphi _1 = \varphi _1(x) $  satisfy
the differential inequality
\begin{eqnarray}
\label{5b}
L  [u ]  + \Delta \varphi _0+ t\Delta \varphi _1 \leq 0 \,,  \quad \mbox{ for all} \quad t \leq T\,,
\end{eqnarray}
and $u$ takes the initial values
\begin{equation}
\label{IC}
u(x,0)=\varphi _0(x)\,, \,\,  u_t(x,0) =\varphi _1(x)   \quad  \mbox{ for all}  \quad  x \in D_0 \subseteq {\mathbb R}^3.
\end{equation}
Then
$
u (x,t) \leq \varphi _0(x)+t \varphi _1(x)\,
$
in the domain of dependence of $D_0$, where $t\leq T $.
\end{theorem}
\medskip

\noindent
{\bf Proof.} For $w=u- \varphi _0-t \varphi _1(x)$ we have
\[
L  [w] = L[u] -L[\varphi _0] - L[t\varphi _1]= L[u]+ \Delta \varphi _0+ t\Delta \varphi _1\leq 0 \,,  \quad \mbox{\rm for all} \quad t \leq T\,,
\]
and the initial conditions
\[
w(x,0)=0\,, \quad  w_t(x,0) =   0   \quad \mbox{\rm for all} \quad x \in D_0 \subseteq {\mathbb R}^3\,.
\]
Then we apply Theorem~1~\cite{Sather}.
\hfill $\square$
\begin{corollary}
Assume that the function $u\in C^2$ satisfies
\[
L  [u ] + \Delta \varphi _0 \leq 0 \,,  \quad \mbox{ for all} \quad t \leq T\,,
\]
and $u$ takes the initial values  (\ref{IC}), where $\varphi _1(x) \leq 0 $ in $D_0$.  Then
$
u (x,t) \leq \varphi _0(x) \,
$
in the domain of dependence of $D_0$, where $t\leq T $.
\end{corollary}
The definition of the superharmonic functions will be used in the next corollaries can be found in \cite{Weinstein}. We are not going to prove the next statements for the less smooth superharmonic functions or for superharmonic function of higher order.
\begin{corollary}
\label{C1.3}
Assume that the function $u$ satisfies
\begin{eqnarray*}
L  [u ]  \leq 0 \,,  \quad \mbox{ for all} \quad t \leq T\,,
\end{eqnarray*}
and $u$ takes the initial values  (\ref{IC}), where $\varphi _1(x) \leq 0 $ in $D_0$.
Suppose that $ \varphi _0  \in C^2$ is  superharmonic in  $D_0 \subseteq {\mathbb R}^3$. Then
$
u (x,t) \leq \varphi _0(x) \,
$
in the domain of dependence of $D_0$, where $t\leq T $.
\end{corollary}

\begin{corollary}
\label{C1.2}
Assume that the function $u$ satisfies
\begin{eqnarray*}
L  [u ]  \leq 0 \,,  \quad \mbox{ for all} \quad t \leq T\,,
\end{eqnarray*}
and $u$ takes the initial values (\ref{IC}).
Suppose that $ \varphi _0 , \varphi _1\in C^2$ are  superharmonic in  $D_0 \subseteq {\mathbb R}^3$. Then
$
u (x,t) \leq \varphi _0(x)+t \varphi _1(x)\,
$
in the domain of dependence of $D_0$, where $t\leq T $.
\end{corollary}
\begin{remark}
The analogous  statements are valid with the subharmonic 
functions $\varphi _0 $, $\varphi _1 \in C^2$.
\end{remark}
We also note that the conditions on the first and second initial data of the solution to the partial differential inequalities (\ref{2}) and (\ref{5b})  are asymmetric. The asymmetry exists also in the Cauchy problem but it reveals itself only in the loss of regularity in one derivative in the Sobolev spaces.
\medskip

Thus, Theorem~\ref{T1.1b}, in particular, gives  sufficient conditions for the solution of  the linear equation to be sign-preserving. If we turn to the linear Klein-Gordon  equation in the Minkowski space
\[
 u_{tt} -  \Delta u  +m^2 u  =   f ,
\]
with $m>0$, then the functional $F(t):=\int_{{\mathbb R}^3} u(x,t)\, dx $ solves the differential equation $F'{}'+m^2 F=\int_{{\mathbb R}^3}f(x,t ) \,dx $.  
The solution $u=u(x,t)$  cannot preserve the sign, for instance, if $ u(x,0)=0$, $f=f(x)$,   
and 
\[
2\left|  \int_{{\mathbb R}^3}f(x,t ) \,dx  \right| <  \left|   \int_{{\mathbb R}^3}u_t(x,0) \,dx  \right|,
\]
since
\[
F(t)
  =
 \frac{1}{m} \sin (mt)   \int_{{\mathbb R}^3}u_t(x,0) \,dx
+ \int_0^t \frac{1}{m}   \sin (m(t- \tau )) \int_{{\mathbb R}^3}f(x,\tau ) \,dx\,d \tau  .
\]
On the other hand, for the the linear Klein-Gordon  operator with the imaginary mass $ L_{KGM}:=\partial_{t }^2 -  \Delta   -M^2$, if 
\begin{equation}
\label{6b}
  L_{KGM}[u] =  f \,,
\end{equation}
then for the functional $F$  with $F(0)=0$, we have
\[
F(t)
  =
\frac{1}{M} \sinh (Mt)  \int_{{\mathbb R}^3}u_t(x,0) \,dx
+ \int_0^t \frac{1}{M}  \sinh (M(t-\tau)) \int_{{\mathbb R}^3}f(x,\tau ) \,dx\,d \tau .
\]
Although the functional $F $ for $f \leq  0$ and $\int_{{\mathbb R}^3}u_t(x,0) \,dx \leq 0 $ is non-positive if $t$ is large,
  we cannot  conclude  that the solution  $u$ is sign preserving. 
\medskip

On the other hand, we can apply the integral transform approach (see   \cite{MN} and references therein) and obtain the following result for the  equation (\ref{6b}).
\begin{theorem}
\label{T1.2}
Assume that  the function $u$ satisfies
\[ 
L_{KGM}[u] \leq 0 \,,  \quad \mbox{  for all} \quad t \leq T\,,
\]
and $L_{KGM}[u] \in C^2$ is a superharmonic in $x$ function.
  Suppose that $ u(x,0)$ and $ u_t(x,0)$ are  superharmonic   non-positive functions  in  $D_0 \subseteq {\mathbb R}^3$. Then  
   \begin{eqnarray}
\label{21bb} 
u(x,t)
&  \leq   &
\int_{ 0}^{t}L_{KGM}[u](x,b)
 \frac{1}{M} \sinh (M(t-b)) \, db   \\
&  &
+  \cosh(Mt) u   (x,0 )  
+\frac{1}{M}\sinh(Mt)u _t(x,0)  \quad \mbox{for all} \quad  t \leq T  \nonumber 
\end{eqnarray}
in the domain of dependence of $D_0$. In particular,
\begin{eqnarray*}
u (x,t) \leq 0 \qquad \mbox{for all} \quad  t \leq T  
\end{eqnarray*}
in the domain of dependence of $D_0$.

\end{theorem}
\medskip

\noindent
{\bf Proof.} If we denote $f:= ( \partial_{t }^2 -  \Delta   -M^2)$, $\varphi _0:= u(x,0)$, and $ \varphi _1:= u_t(x,0)$, then according to the integral transform approach formulas \cite{MN} we can write 
 \begin{eqnarray*} 
u(x,t)
& = &
\int_{ 0}^{t} db
  \int_{ 0}^{t-b }   I_0\left(M\sqrt{(t-b)^2-r^2}\right)   v_f(x,r;b ) \, dr \\
&  &
+ v_{\varphi _0}(x,t )+\int_{ 0}^{t}    \frac{\partial }{\partial t}  I_0\left(M\sqrt{t^2-r^2}\right)        v_{\varphi _0}(x,r  ) \,  dr \\
&  &
+\int_{ 0}^{t}   I_0\left(M\sqrt{t^2-r^2}\right)   v_{\varphi _1}(x,r  )\,  dr, \,\, x \in {\mathbb R} , \,\, t> 0,
\end{eqnarray*}
where the function
$v(x,t;b)$
is the solution to the Cauchy problem for the  wave equation
\[
v_{tt} -   \bigtriangleup v  =  0 \,, \quad v(x,0;b)=f(x,b)\,, \quad v_t(x,0)= 0\,,
\]
while   $v_{\varphi } $
is the solution  of the Cauchy problem
\[
v_{tt} -   \bigtriangleup v  =  0 \,, \quad v(x,0)=\varphi  (x)\,, \quad v_t(x,0)= 0\,.
\]
Here  $I_0(z)$ is the modified Bessel function of the first kind. Then the statement of this theorem follows from Teorem~\ref{T1.1b} and the properties of the function $I_0\left(z\right) $. Indeed, due to Corollary~\ref{C1.3}, we have
\[
v_f(x,r;b ) \leq f(x,b)\,, \quad v_{\varphi _0}(x,r  )\leq \varphi _0 (x)\,, \quad v_{\varphi _1}(x,r  )\leq \varphi _1 (x)
\]
for all corresponding $x,r$, and $b$. The function  $ I_0(z)$ is positive while $ I_0'(z)$ is non-negative for $z>0$.  Thus,
the inequality 
 \begin{eqnarray*} 
u(x,t)
& \leq  &
\int_{ 0}^{t} db \, f(x,b)
  \int_{ 0}^{t-b }   I_0\left(M\sqrt{(t-b)^2-r^2}\right)     \, dr \\
&  &
+  \varphi _0 (x  )+\varphi _0 (x  )\int_{ 0}^{t}    \frac{\partial }{\partial t}  I_0\left(M\sqrt{t^2-r^2}\right)          \,  dr \\
&  &
+ \varphi _1 (x )\int_{ 0}^{t}   I_0\left(M\sqrt{t^2-r^2}\right)  \,  dr, \,\, x \in {\mathbb R} , \,\, t> 0,
\end{eqnarray*}
and the result of the integrations prove  theorem. 
\hfill $\square$

 On the other hand, in order to  prove a sign changing property of the solutions to the semilinear equations for those no explicit formulas are available, the $F$-functional method can be applied.  For details see \cite{CPDE2012}.

\section{The maximum principle in the de~Sitter space-time}
\label{S3}

For the hyperbolic equation with variable coefficients the maximum principle is known only in the one dimensional case (see, e.g., \cite{Protter})
and for Euler-Poisson-Darboux equation \cite{Weinstein}.
We consider the linear part of the equation
\begin{equation}
\label{K_G_Higgs}
u_{tt} - e^{-2t} \bigtriangleup u  - M^2   u=  - e^{\frac{n}{2}t}V'(e^{-\frac{n}{2}t}u ),
\end{equation}
with $M\geq 0 $ and the potential function $V=V(\psi ) $. If we denote  the non-covariant Klein-Gordon operator in the de~Sitter space-time
\[
L_{KGdS}:= \partial_t^2 - e^{-2t} \bigtriangleup    - M^2     ,
\]
then (\ref{K_G_Higgs}) can be written as follows:
\[
L_{KGdS}[u]=  - e^{\frac{n}{2}t}V'(e^{-\frac{n}{2}t}u )\,.
\]
The  equation (\ref{K_G_Higgs}) covers two important cases. The first one is the Higgs boson equation (\ref{HBE}) that leads to
(\ref{K_G_Higgs}) if  $\psi = e^{-\frac{n}{2}t} u$.  Here $V'(\psi )=\lambda \psi ^3 $
and $M^2=  \mu ^2+ n^2/4 $ with $\lambda >0 $ and $\mu >0 $, while $n=3$. The second case is the case of the covariant Klein-Gordon equation
\[
\psi _{tt} + n\psi _t- e^{-2t} \bigtriangleup \psi   + m^2  \psi =  -  V'(\psi  ),
\]
with small physical  mass, that is $0 \leq m  \le  {n }/{2}$. For the last case
$ M^2=  {n^2}/{4}-m^2$. It is evident that the last equation is related to the equation (\ref{K_G_Higgs}) via transform $\psi =e^{-\frac{n}{2}t}u $.
\smallskip

It is known that the Klein-Gordon quantum fields   whose squared physical masses are negative
(imaginary mass)  represent tachyons. (See, e.g., \cite{B-F-K-L}.)
   In \cite{B-F-K-L}   the Klein-Gordon equation with imaginary mass is considered. It is shown that localized disturbances spread with at most the speed of light, but grow exponentially. The conclusion is made that    free tachyons have to be rejected   on stability grounds.
\medskip

The Klein-Gordon quantum fields on the de~Sitter manifold  with
imaginary mass 
present   scalar tachyonic  quantum fields.    Epstein and Moschella \cite{Epstein-Moschella} give an exhaustive   study of   scalar tachyonic quantum
fields which are linear Klein-Gordon quantum
fields on the de Sitter manifold whose  masses   take an infinite
set of discrete values  
$m^2=-k(k+n)$, $k=0,1,2,\ldots$. The corresponding linear equation is
\begin{eqnarray*}
&  &
\psi _{tt} +   n  \psi  _t - e^{-2 t} \Delta  \psi   +  m^2  \psi  =  0\,.
 \end{eqnarray*}
If $n$ is an odd number, then
$m$ takes value at  the  knot points set  \cite{JMP2013}.

The nonexistence of a global in time  solution of the semilinear Klein-Gordon massive tachyonic (self-interacting quantum
fields)  equation in the de~Sitter space-time is proved in \cite{yagdjian_DCDS}.  
More precisely, consider the   semilinear equation
\begin{eqnarray*}
&  &
\psi _{tt} +   n  \psi  _t - e^{-2 t} \Delta  \psi   -  m^2  \psi  =  c|\psi|^{1+\alpha} \,,
 \end{eqnarray*}
which is   commonly used model for general nonlinear problems. Then, according to Theorem~1.1~\cite{yagdjian_DCDS}, if $c\not= 0$, $\alpha >0 $,
and $m \not= 0$, then for every positive numbers $\varepsilon  $ and $s$ there exist  functions $\psi _0 $,
$\,  \psi _1  \in C_0^\infty ( {\mathbb R}^n)$ such that $ \|\psi _0 \|_{H_{(s)}( {\mathbb R}^n)} + \|\psi _1 \|_{H_{(s)}( {\mathbb R}^n)} \leq \varepsilon $ but the solution $\psi =\psi (x,t) $ with the  initial values
\[
  \psi (x,0)= \psi _0 (x)\,, \quad \psi _t(x,0)= \psi _1 (x)\,,
\]
blows up in finite time. This implies also blowup for the sign-preserving solutions of the equation
\begin{eqnarray*}
&  &
\psi _{tt} +   n  \psi  _t - e^{-2 t} \Delta  \psi   -  m^2  \psi  =  c|\psi|^{ \alpha} \psi  \,.
 \end{eqnarray*}

\medskip

The next theorem gives certain  kind of maximum principle for the  non-covariant  Klein-Gordon equation in the de~Sitter space-time.
Define the ``forward light cone'' $D^{dS}_+ (x_0,t_0) $
and   the ``backward light cone'' $D^{dS}_- (x_0,t_0) $,  in the de~Sitter space-time for the point $(x_0,t_0 ) \in {\mathbb R}^{n+1}$,
  as follows
\begin{eqnarray*}
D^{dS}_\pm (x_0,t_0)
& := &
\Big\{ (x,t)  \in {\mathbb R}^{n+1}  \, ; \,
|x -x_0 | \leq \pm( e^{-t_0} - e^{-t })
\,\Big\} \,.
\end{eqnarray*}
For the domain $D_0\subseteq {\mathbb R}^n$ define dependence domain of $D_0 $ as follows:
\begin{eqnarray*}
D^{dS}(D_0)
& := &
\displaystyle \bigcup_{x_0 \in {\mathbb R}^n,\,\, t_0 \in [0,\infty) }    \left\{ D^{dS}_- (x_0,t_0) \,;\,  D^{dS}_- (x_0,t_0) \cap \{t=0\} \subset D_0 \right\}\,.
\end{eqnarray*}
\begin{theorem}
\label{T2}
Assume that  $M >1  $ and the function $u$ satisfies
\[ 
L_{KGdS}[u] \leq 0 \,,  \quad \mbox{  for all} \quad t \leq T\,,
\]
and $L_{KGdS}[u] \in C^2$ is a superharmonic in $x$ function.
  Suppose that $ u(x,0)$ and $ u_t(x,0)$ are  superharmonic   non-positive functions in  $D_0 \subseteq {\mathbb R}^3$. 
Then  
\begin{eqnarray}
\label{21bbc}
u(x,t)
&  \leq   &
\int_{ 0}^{t} L_{KGdS}[u](x,b)
 \frac{1}{M} \sinh (M(t-b)) \, db  
+  \cosh(Mt) u   (x,0 )  
+\frac{1}{M}\sinh(Mt)u _t(x,0)   
\end{eqnarray}
for all   $t\in [\ln (M/(M-1)),T]$ in the domain of dependence of $D_0$. In particular,
\begin{eqnarray}
\label{21b}
u (x,t) \leq 0 \qquad \mbox{for all} \quad  t\in [\ln (M/(M-1)),T] 
\end{eqnarray}
in the domain of dependence of $D_0$.

If $u(x,0) \equiv 0 $, then the statements (\ref{21bbc}),(\ref{21b}) hold  also for all $t \in [0,T] $ and  each $M\geq 0 $.
\end{theorem}
\medskip

\noindent
{\bf Proof.} 
We are going to apply the integral transform and  the  kernel functions $E(x,t;x_0,t_0;M) $, $K_0(z,t;M)   $,    and $K_1(z,t;M) $ from \cite{JMAA_2012}.
 First we  introduce the function
\begin{eqnarray*}
E(x,t;x_0,t_0;M)
&  =  &
 4 ^{-M}  e^{ M(t_0+t) } \Big((e^{-t }+e^{-t_0})^2 - (x - x_0)^2\Big)^{-\frac{1}{2}+M    } \\
&  &
\times F\Big(\frac{1}{2}-M   ,\frac{1}{2}-M  ;1;
\frac{ ( e^{-t_0}-e^{-t })^2 -(x- x_0 )^2 }{( e^{-t_0}+e^{-t })^2 -(x- x_0 )^2 } \Big) . \nonumber
\end{eqnarray*}
Here $F\big(a, b;c; \zeta \big) $ is the hypergeometric function. (See, e.g., \cite{B-E}.) Next
we define the kernels  $K_0(z,t;M)   $    and $K_1(z,t;M) $ by
\begin{eqnarray}
\label{K0}
K_0(z,t;M) 
&   :=  &
- \left[  \frac{\partial }{\partial b}   E(z,t;0,b;M) \right]_{b=0}   \\
&  = &
4 ^ {-M}  e^{ t M}\big((1+e^{-t })^2 - z^2\big)^{ -\frac{1}{2}+ M    } \frac{1}{  (1-e^{ -t} )^2 -  z^2  } \nonumber \\
&   &
\times  \Bigg[  \big(  e^{-t} -1 +M(e^{ -2t} -      1 -  z^2) \big)
F \Big(\frac{1}{2}-M   ,\frac{1}{2}-M  ;1; \frac{ ( 1-e^{-t })^2 -z^2 }{( 1+e^{-t })^2 -z^2 }\Big)  \nonumber \\
&  &
+   \big( 1-e^{-2 t}+  z^2 \big)\Big( \frac{1}{2}+M\Big)
F \Big(-\frac{1}{2}-M   ,\frac{1}{2}-M  ;1; \frac{ ( 1-e^{-t })^2 -z^2 }{( 1+e^{-t })^2 -z^2 }\Big) \Bigg] \nonumber
\end{eqnarray}
and $K_1(z,t;M)   :=
  E(z ,t;0,0;M) $, that is,
\begin{eqnarray*}
K_1(z,t;M)
& = &
  4 ^{-M} e^{ Mt }  \big((1+e^{-t })^2 -   z  ^2\big)^{-\frac{1}{2}+M    } \\
&  &
\times F\left(\frac{1}{2}-M   ,\frac{1}{2}-M  ;1;
\frac{ ( 1-e^{-t })^2 -z^2 }{( 1+e^{-t })^2 -z^2 } \right), \, 0\leq z\leq  1-e^{-t},
 \end{eqnarray*}
respectively. These kernels have  been introduced and used in \cite{Yag_Galst_CMP,yagdjian_DCDS} in the  representation of the solutions of the Cauchy problem. The positivity of the kernels $E $, $K_0 $, and $K_1 $ is proved in the next section.
\smallskip

The solution $u= u(x,t)$ to the Cauchy problem
\[
u_{tt} - e^{-2t}\Delta u -M^2 u= f ,\quad u(x,0)= 0  , \quad u_t(x,0)=0,
\]
with \, $ f \in C^\infty ({\mathbb R}^{n+1})$\, and with   vanishing
initial data is given  \cite{JMAA_2012} by the next expression
\begin{eqnarray}
\label{1.29small}
u(x,t)
&  =  &
2   \int_{ 0}^{t} db
  \int_{ 0}^{ e^{-b}- e^{-t}} dr  \,  v(x,r ;b) E(r,t; 0,b;M)  ,
\end{eqnarray}
where the function
$v(x,t;b)$
is a solution to the Cauchy problem for the  wave equation
\begin{equation}
\label{1.6c}
v_{tt} -   \bigtriangleup v  =  0 \,, \quad v(x,0;b)=f(x,b)\,, \quad v_t(x,0)= 0\,.
\end{equation}
If the  superharmonic function $f$ is also non-positive, $f(x,t) \leq 0 $, then    due to Corollary~\ref{C1.2} we conclude
\begin{eqnarray*}
 v(x,r ;b)
&   \leq   &
f(x,b) \leq 0,
\end{eqnarray*}
in the domain of dependence of $D_0$. It follows
\begin{eqnarray*}
   \int_{ 0}^{t} db
  \int_{ 0}^{ e^{-b}- e^{-t}} dr  \,  v(x,r ;b) E(r,t; 0,b;M)  
  & \leq & 
  \int_{ 0}^{t}  f(x,b) \, db
  \int_{ 0}^{ e^{-b}- e^{-t}}  E(r,t; 0,b;M) \,dr   \\ 
  & \leq & 
 \int_{ 0}^{t} f(x,b)
 \frac{1}{2M} \sinh (M(t-b)) \, db    \leq 0\,,
\end{eqnarray*}
provided that $E(r,t; 0,b;M)\geq 0 $.

The solution $u=u (x,t)$ to the Cauchy problem
\begin{equation}
\label{CPu}
u_{tt}-  e^{-2t} \bigtriangleup u -M^2 u =0\,, \quad u(x,0)= u_0 (x)\, , \quad u_t(x,0)=u_1 (x)\,,
\end{equation}
with \, $u_0 $,  $ u_1 \in C_0^\infty ({\mathbb R}^n) $, $n\geq 2$, can be represented \cite{JMAA_2012}  as follows:
\begin{eqnarray*}
u(x,t)
& = &
 e ^{\frac{t}{2}} v_{u_0}  (x, \phi (t))
+ \, 2\int_{ 0}^{1} v_{u_0}  (x, \phi (t)s) K_0(\phi (t)s,t;M)\phi (t)\,  ds  \nonumber \\
& &
+\, 2\int_{0}^1   v_{u_1 } (x, \phi (t) s)
  K_1(\phi (t)s,t;M) \phi (t)\, ds
, \quad x \in {\mathbb R}^n, \,\, t>0\,,
\end{eqnarray*}
where $\phi (t):=  1-e^{-t} $.
 Here, for $\varphi \in C_0^\infty ({\mathbb R}^n)$ and for $x \in {\mathbb R}^n$,
the function $v_\varphi  (x, \phi (t) s)$  coincides with the value $v(x, \phi (t) s) $
of the solution $v(x,t)$ of the Cauchy problem
\[
v_{tt} -   \bigtriangleup v  =  0 \,, \quad v(x,0)=\varphi (x)\,, \quad v_t(x,0)= 0\,.
\]

For the function  $u _1   $, which is  superharmonic, from Corollary~\ref{C1.3} we conclude
\[
v_{u _1} (x,r) \leq  u _1(x)  \,.
\]
It follows
\begin{eqnarray*}
 2\int_{0}^{\phi (t)}   v_{u _1} (x,  r)
  K_1(r,t;M)  \, dr 
& \leq &
 2 u _1(x) \int_{0}^{\phi (t)}
  K_1(r,t;M)  \, dr \\
& = &
\frac{1}{M}\sinh(Mt)u _1(x)\,.
\end{eqnarray*}
since
$
K_1(r,t;M) \geq 0$.  
In particular, if $u _1(x) \leq 0$, then
\begin{eqnarray*}
&    &
 2\int_{0}^{\phi (t)}   v_{u _1 } (x,  r)
  K_1(r,t;M)  \, dr \leq 0\,.
\end{eqnarray*}

Further, if $u_0 \in C^2$ is superharmonic, that is $\Delta u_0 \leq 0 $, then, according to  Corollary~\ref{C1.3},
$
(\partial_t^2-\Delta ) v_{u_0} = 0
$
implies $v_{u_0}(x,t) \leq u_0(x)$. Consequently, if $ M>1$,  then $ K_0(r,t;M) \geq 0$ for all $t \in [\ln(M/(M-1)),T]$,
and 
\begin{eqnarray*}
&   &
  e ^{\frac{t}{2}} v_{u_0}  (x, \phi (t))
+ \, 2\int_{ 0}^{\phi (t)} v_{\varphi_0}  (x, r) K_0(r,t;M) \,  dr   \\
& \leq   &
  u_0  (x ) \left[  e ^{\frac{t}{2}}
+ \, 2\int_{ 0}^{\phi (t)}   K_0(r,t;M) \,  dr \right]= \cosh(Mt) u_0  (x )\,.
\end{eqnarray*}
 Theorem~\ref{T2} is proved. \hfill $\square$
\medskip

We do not know if the condition  of superharmonicity can be relaxed.

\section{The positivity of the kernel functions $E$, $K_0$ and $K_1$}
\label{S5}

\begin{proposition}
\label{P1}
Assume that $M\geq  0$. 
Then
\begin{eqnarray*}
&    &
 E(r,t; 0,b;M) \geq 0 ,  \quad \mbox{  for all}\quad   0  \leq b \leq t,\quad r\leq e^{-b}-e^{-t}\,,\quad  t \in [0,\infty)\,, \\
& &
K_1(r,t;M)  \geq 0\quad \mbox{  for all}\quad  r\leq 1-e^{-t}, \quad  t \in [0,\infty)\,.
\end{eqnarray*}
If we assume that $M>1$, then
\begin{eqnarray*}
&   &
K_0(r,t;M) \geq 0   \quad \mbox{\rm for all}\,\, r\leq 1-e^{-t} \,\, \mbox{ and  for all }\,\, t >   \ln \frac{M}{M-1}      \,.
\end{eqnarray*}
 \end{proposition}
\medskip

\noindent
{\bf Proof.} Indeed, for $ 0  \leq b \leq t$ and  $r\leq e^{-b}-e^{-t}$ we have
\begin{eqnarray*}
E(r,t;0,b;M)
&  =  &
 4 ^{-M}  e^{ M(b+t) } \Big((e^{-t }+e^{-b})^2 - r^2\Big)^{-\frac{1}{2}+M    } \\
&  &
\times F\Big(\frac{1}{2}-M   ,\frac{1}{2}-M  ;1;
\frac{ ( e^{-b}-e^{-t })^2 -r^2 }{( e^{-b}+e^{-t })^2 -r^2 } \Big)  .
\end{eqnarray*}
For   $M \geq 0$   the parameters $a=b= {1}/{2}-M$ and  $c=1$ $ $  of the   
function $F(a,b;c;z)$ satisfy the relation $a+b\leq c$. Then, we denote
\[
z    :=  
 \frac{ ( e^{-b}-e^{-t })^2 -r^2 }{( e^{-b}+e^{-t })^2 -r^2 } \leq 1 , \quad 0  \leq b \leq t,\quad r\leq e^{-b}-e^{-t} \,.
\]
Hence, it remains to check the sign of the  
function $F ( a   ,a  ;1; z  ) $ with parameter $a \leq 1/2  $ and $ z \in(0,1)$. If $a$ is not a non-positive integer then the series
\begin{eqnarray*}
&  &
F ( a   ,a  ;1; x  )= \sum_{n=0}^\infty \frac{[(a)_n]^2}{[n!]^2}x^n \,,\quad (a)_n:= a(a+1)\cdots(a+n-1)  \,,
\end{eqnarray*}
is a convergent series for all $x \in [0,1)$. If $a$ is negative integer, $a=-k$, then $F (a,a  ;1;x ) $ is polynomial with the positive coefficients:
\begin{eqnarray*}
&  &
F ( a   ,a  ;1; x  )= \sum_{n=0}^k \frac{[(a)_n]^2}{[n!]^2}x^n \,.
\end{eqnarray*}
Since $K_1(z,t;M)   :=   E(z ,t;0,0;M) $, the first two statements of the proposition are proved.

In order to verify the last statement
it suffices to verify the inequality $K_0(r,t;M) >0 $,
where $ r \in (0,1)$. Denote    $M=(2k+1)/2$. Then $ {1}/{2}-M= -k <0$ and 
we can write (\ref{K0}) in the equivalent form  as follows
\begin{eqnarray*}
K_0(z,t;M)
&  = &
- \Bigg[  \frac{\partial }{\partial b} \Bigg\{   4 ^{-M}  e^{ M(b+t) } \Big((e^{-t }+e^{-b})^2 - r^2\Big)^{-\frac{1}{2}+M    }  \Bigg\}\Bigg]_{b=0} \nonumber \\
&  &
\times F\Big(\frac{1}{2}-M   ,\frac{1}{2}-M  ;1;
\frac{ ( e^{-b}-e^{-t })^2 -r^2 }{( e^{-b}+e^{-t })^2 -r^2 } \Big)  \nonumber \\
&   &
-   4 ^{-M}  e^{ M(b+t) } \Big((e^{-t }+e^{-b})^2 - r^2\Big)^{-\frac{1}{2}+M    }  \nonumber \\
&  &
\times \Bigg[  \frac{\partial }{\partial b} \Bigg\{ F\Big(\frac{1}{2}-M   ,\frac{1}{2}-M  ;1;
\frac{ ( e^{-b}-e^{-t })^2 -r^2 }{( e^{-b}+e^{-t })^2 -r^2 } \Big) \Bigg\}\Bigg]_{b=0} \,. \nonumber
\end{eqnarray*}
Then we use the relation (20)~\cite[Sec.2.8]{B-E}:
\begin{eqnarray*}
K_0(z,t;M)
&  = &
- 4^{-(k+\frac{1}{2})} e^{(k+\frac{1}{2}) t} \Big((e^{-t }+1)^2 - r^2\Big)^{k} \nonumber \\
&  &
\times \Bigg\{  \frac{\left(-e^{2 t} \left((k+\frac{1}{2}) r^2+(k+\frac{1}{2})-1\right)+(k+\frac{1}{2})+e^t\right)}{\left(r^2-1\right) e^{2 t}-2 e^t-1}\nonumber \\
&  &
\times F\Big(-k  ,-k ;1;
\frac{ ( 1-e^{-t })^2 -r^2 }{( 1+e^{-t })^2 -r^2 } \Big)  \nonumber \\
&   &
-       k^2  F\Big(1-k ,1-k ;2;
\frac{ ( 1-e^{-t })^2 -r^2 }{( 1+e^{-t })^2 -r^2 } \Big) \Bigg\}\,.   \nonumber
\end{eqnarray*}
Thus
\begin{eqnarray*}
&  &
K_0 \left(r,t;\frac{1}{2}+k \right) \\
& = &
 4^{-k-1} e^{\left(k+\frac{1}{2}\right) t} \big((1+e^{-t })^2 - r^2\big)^{  k -2  }   \\
&  &
\times  \Bigg[8 k^2 e^t \left(\left(r^2+1\right) e^{2 t}-1\right)
F \left(1-k,1-k;2;\frac{ ( 1-e^{-t })^2 -r^2 }{( 1+e^{-t })^2 -r^2 }\right)\\
&  &
+ \left[(1+e^t)^2-r^2 \right] \left(e^{2 t} \left(2 k \left(r^2+1\right)+r^2-1\right)-2 k-2 e^t-1\right) \\
&  &
\times F \left(-k,-k;1;\frac{ ( 1-e^{-t })^2 -r^2 }{( 1+e^{-t })^2 -r^2 }\right)\Bigg]  \,.
\end{eqnarray*}
Consider the factor
\begin{eqnarray}
\label{20}
&  &
\Bigg[8 k^2 e^t \left(\left(r^2+1\right) e^{2 t}-1\right)
F \left(1-k,1-k;2;\frac{ ( 1-e^{-t })^2 -r^2 }{( 1+e^{-t })^2 -r^2 }\right)\\
&  &
+ \left[(1+e^t)^2-r^2 \right] \left(e^{2 t} \left(2 k \left(r^2+1\right)+r^2-1\right)-2 k-2 e^t-1\right) \nonumber  \\
&  &
\times F \left(-k,-k;1;\frac{ ( 1-e^{-t })^2 -r^2 }{( 1+e^{-t })^2 -r^2 }\right)\Bigg] \,.\nonumber
\end{eqnarray}
The functions $F \left( -k, -k;1;z\right) $ and $F \left( 1-k,1-k;2;z\right) $ are defined as follows
\begin{eqnarray*}
F \left( -k, -k;1;z\right)
& = &
 \sum_{n=0}^\infty  \frac{[ (-k)(-k+1)\cdots(-k+n-1)]^2}{[n!]^2}z^n
 \,,\\
F \left(1-k,1-k;2;z\right)
& = &
\sum_{n=0}^\infty  \frac{ [ (1-k)_n ]^2}{ [n!]^2 (n+1) }z^n \,.
  \end{eqnarray*}
  Here we have denoted
  \[
  z:= \frac{ ( 1-e^{-t })^2 -r^2 }{( 1+e^{-t })^2 -r^2 } \in [0,1] \quad \mbox{\rm for all} \,\,t \in [0,\infty), \,\, r \in (0,1-e^{-t })\,.
  \]
  Thus,
 \begin{eqnarray*}
F \left( -k, -k;1;z\right)
& \geq  &
 1 \quad \mbox{\rm for all} \,\, z \in [0,1)\,,\\
F \left(1-k,1-k;2;z\right)
& \geq  &
 1 \quad \mbox{\rm for all} \,\, z \in [0,1)\,.
  \end{eqnarray*}
On the other hand,
\begin{eqnarray*}
&  &
8 k^2 e^t \left(\left(r^2+1\right) e^{2 t}-1\right)>1  \quad \mbox{\rm for all} \,\, r \in [0,1] \quad \mbox{\rm and } \,\,t \geq \frac{M}{M-1}  \,.
\end{eqnarray*}
Next we   check the sign of the function
\begin{eqnarray*}
&  &
\left[(1+e^t)^2-r^2 \right]
\left(e^{2 t} \left(2 k \left(r^2+1\right)+r^2-1\right)-2 k-2 e^t-1\right)\,.
\end{eqnarray*}
Since $ \left[(1+e^t)^2-r^2 \right] \geq 3$, we consider the second factor only. We set    $x:=e^t>1$ and $y:=r^2  \in [0,1]$; then we have the polynomial
\begin{eqnarray*}
P(x,y)
& = &
 x^{2} \left(2 k \left(y+1\right)+y-1\right)-2 k-2 x-1   \,.
\end{eqnarray*}
It follows
\begin{eqnarray*}
\partial_y P(x,y)
& = &
2x^2M>0 \,.
\end{eqnarray*}
On the other hand,   
\begin{equation}
\label{32}
P(x,0)
  =   2[ x^{2} (  M  -1)- x - M] >0\,,
\end{equation}
and for $M=1$ the last inequality false since $x>0$. For $M>1$ the inequality (\ref{32}) holds whenever $x> M/(M-1)$. 
It follows
\begin{eqnarray*}
P(e^t,r^2)
& > &
const >0\quad  \forall \, r^2 \in [0,1] \quad \mbox{\rm and  } \forall\,t \in (\ln (M/(M-1)),\infty)\,.
\end{eqnarray*}
Since all terms of (\ref{20}) are positive,  the proposition is proved. \hfill $\square$
\medskip

\begin{remark}
The graph of the $K_0(r,t; \frac{3}{4})$ shows that the  $K_0$  changes a sign.
\begin{center}
\begin{figure}[h]
\label{Fig1}
\includegraphics[width=0.32\textwidth]{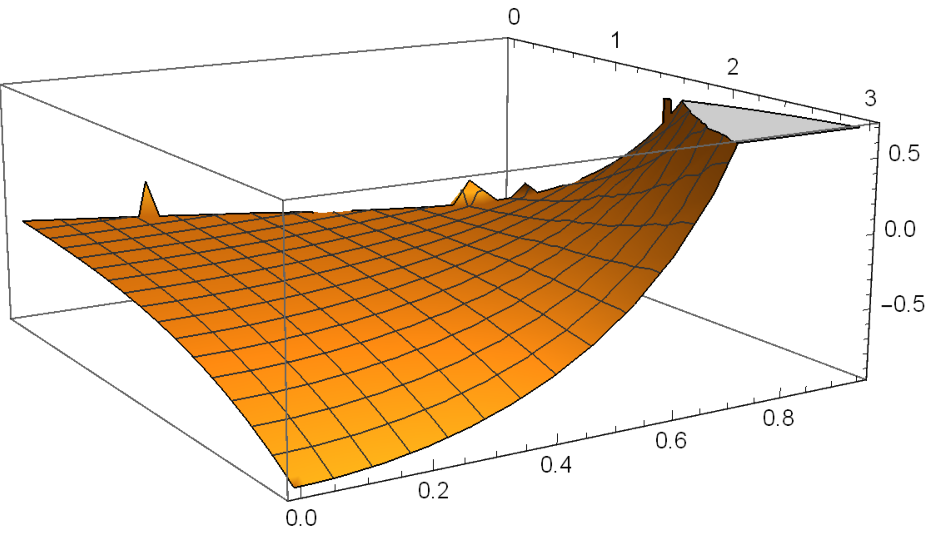}
\includegraphics[width=0.32\textwidth]{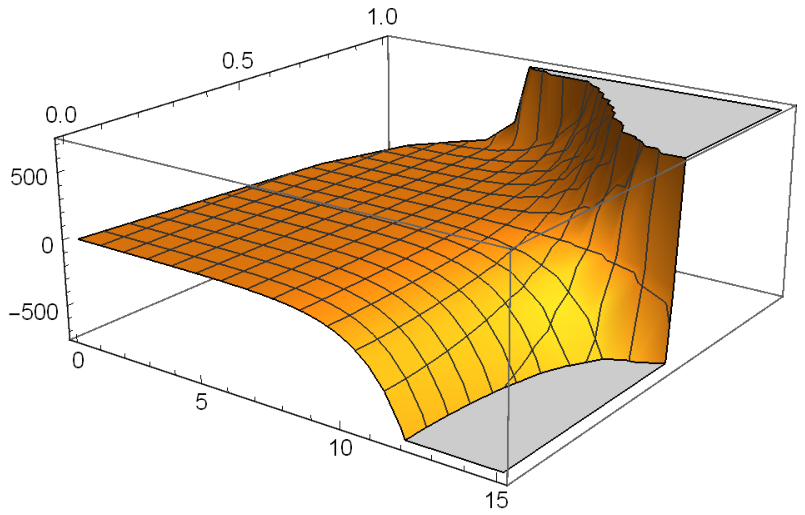}
\includegraphics[width=0.32\textwidth]{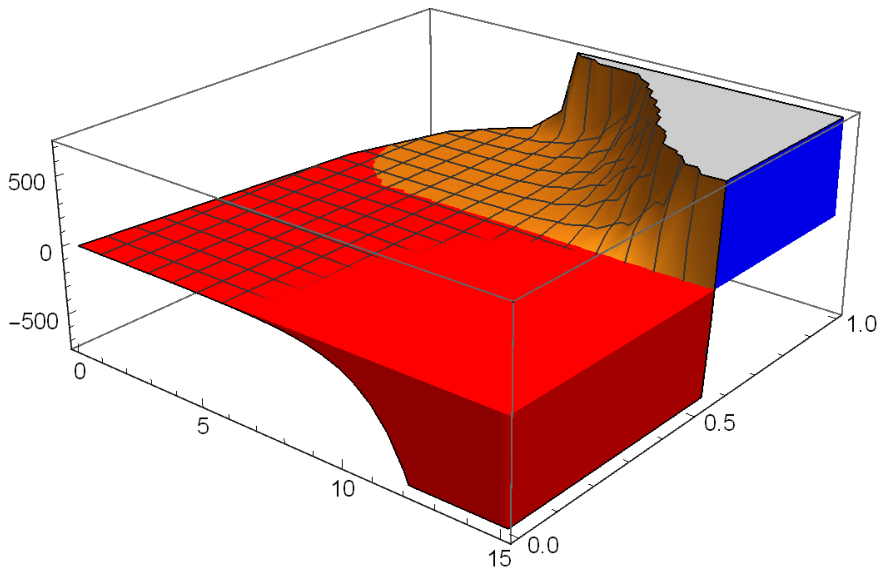}\\
\caption{\small The graph of $K_0 \left(z,t,\frac{3}{4}\right)$, $\,t\in (0,3)$ and $\,t\in (0,15)$, $\,z\in (0,1-\exp (-t))$}
\end{figure}
\end{center}
\end{remark}
\begin{remark}
The graph of the $K_0(r,t; \frac{1}{6})$ shows that the  $K_0$ does not  change a sign.
\begin{center}
\begin{figure}[h]
\label{Fig1b}
\hspace{1.8cm}\includegraphics[width=0.32\textwidth]{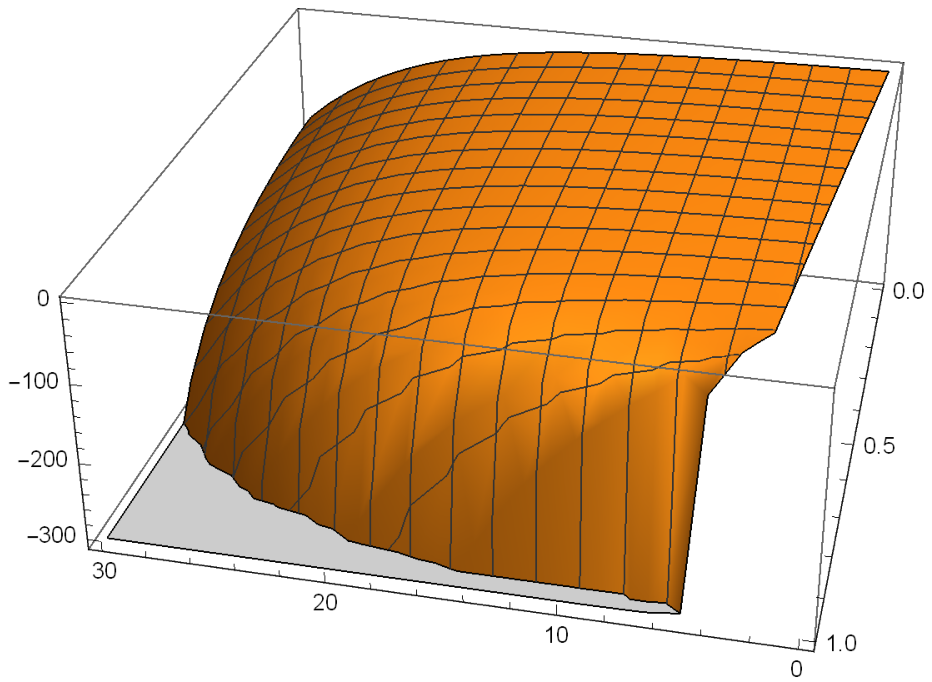}\hspace{2.9cm}
\includegraphics[width=0.32\textwidth]{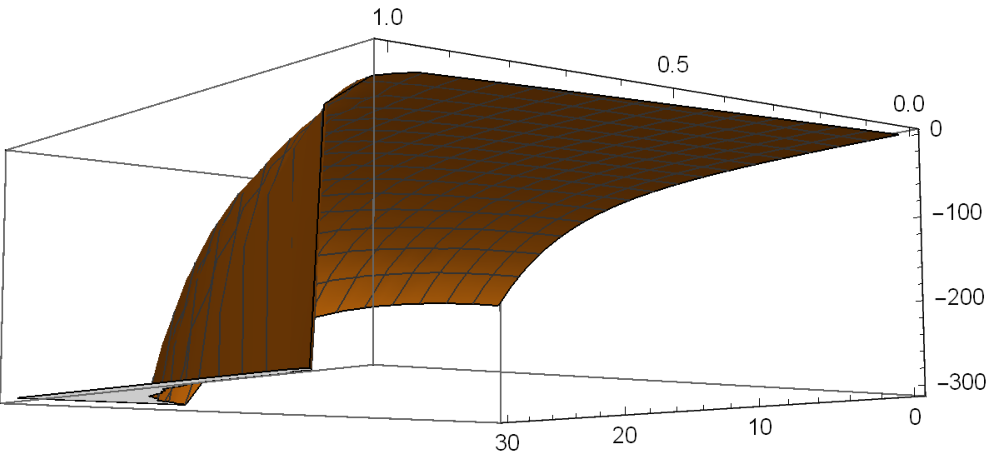}\\
\caption{\small The graph of $K_0 \left(z,t,\frac{1}{6}\right)$, $\,t\in (0,30)$, $\,z\in (0,1-\exp (-t))$}
\end{figure}
\end{center}
\end{remark}

For $M=1/2$ the kernels are (see \cite{JMP2013})
\[
 E \left( r,t;0,b;\frac{1}{2} \right)= \frac{1}{2}e^{\frac{1}{2}(b+t)},\quad  K_0\left(r,t;\frac{1}{2}\right)= -\frac{1}{4}e^{\frac{1}{2} t } ,\quad
K_1\left(r,t;\frac{1}{2}\right)=  \frac{1}{2}e^{\frac{1}{2} t } \,
\]
and the solution can be written as follows
\begin{eqnarray*} 
u(x,t) 
&  =  &
 \int_{ 0}^{t} db
  \int_{ 0}^{ e^{-b}- e^{-t}} dr  \,  e^{\frac{1}{2}(b+t)} v(x,r ;b)
+ e ^{\frac{1}{2}t} v_{u_0}  (x, \phi (t)) \\
&  &
  -\,  \frac{1}{2} e^{\frac{1}{2} t }\int_{ 0}^{\phi (t)} v_{u_0}  (x, r) \,  dr  
+\,   e^{\frac{1}{2} t }\int_{0}^{\phi (t)}    v_{u_1 } (x, r)
\, dr,  
\end{eqnarray*}
where $x \in {\mathbb R}^n$, $t>0$. In particular, if $ f=0$ and $u_1=0 $, then 
\[
u(x,t)
   =    e ^{\frac{1}{2}t} v_{u_0}  (x, \phi (t))
  -\,  \frac{1}{2} e^{\frac{1}{2} t }\int_{ 0}^{\phi (t)} v_{u_0}  (x, r) \,  dr,  \quad x \in {\mathbb R}^n, \,\, t>0,
\]
 solves   (\ref{CPu}).  The second term of the last expression is the so-called tail. The tail is of considerable interest in many aspects in physics, and, in particular, in the General Relativity \cite{Sonego}.
\begin{remark}
If we assume that
$ u(x,t) \leq 0$, 
then
\begin{eqnarray*}
v_{u_0}  (x, s)
  -\,  \frac{1}{2}\int_{ 0}^{s} v_{u_0}  (x, r) \,  dr   \leq 0\quad \mbox{ for all} \quad  s \in[0,1]\,,
\end{eqnarray*}
and the Gronwall's lemma implies
$
v_{u_0}  (x, s) \leq 0\,.
$
\end{remark}
The converse statement is not true in general. Indeed, according to the Figure~3, for  $v(s) = -e^{-\frac{s^2}{1.2\, -s^3}} <0 $, the   function
\begin{eqnarray*}
v   (  s)
  -\,  \frac{1}{2}\int_{ 0}^{s} v   (  r) \,  dr    \,,
\end{eqnarray*}
is positive when $s \geq 0.9 $:
\begin{center}
\begin{figure}[h]
\label{F3}
\hspace{2.8cm}\includegraphics[width=0.35\textwidth]{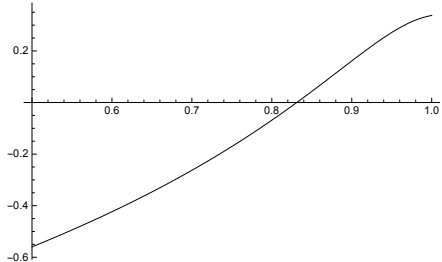}
\caption{\small The graph of $v   (  s)   -\,  \frac{1}{2}\int_{ 0}^{s} v   (  r) \,  dr   $ with $ v(s) = -e^{-\frac{s^2}{1.2\, -s^3}} $}
\end{figure}
\end{center}
If $u_0=u_0(x) $ is harmonic function in ${\mathbb R}^n $, then $v_{u_0}  (x, r) = u_0 (x )$ and
\begin{eqnarray*}
u(x,t)
&  =  &
\cosh \left(\frac{1}{2}t \right)u_0(x)
, \quad x \in {\mathbb R}^n, \,\, t>0\,.
\end{eqnarray*}

\begin{remark}
We do not know if the value $M=1$, that is $m=\sqrt{5}/2$, has some physical significance similar to one when $M=1/2$, that is $m=\sqrt{2} $,  which is the end point of the Higuchi bound \cite{JMAA_2012}.
\end{remark}

\begin{conjecture}
Assume that $M \in[0, 1/2]$.
Then
\begin{eqnarray*}
&   &
K_0(r,t;M) \leq  0   \quad \mbox{\rm for all}\,\, r\leq 1-e^{-t}  \quad \mbox{ and  for all }\,\, t>0  \,.
\end{eqnarray*}
\end{conjecture}
\medskip

\section{The sign-changing solutions for the semilinear Klein-Gordon  equation in the de~Sitter space  with Higgs potential}
\label{S6}

We are interested in sign-changing solutions of the  equation for the Higgs real-valued scalar    field  in the
  de~Sitter space-time
\begin{equation}
\label{1.1}
\label{Higgs_eq}
  \psi_{tt} +   3   \psi_t - e^{-2 t} \Delta  \psi =    \mu^2  \psi  -\lambda   \psi^ 3 \,.
\end{equation}
The constants $\psi =  \pm  {\mu }/{\sqrt{\lambda }} $ are non-trivial  real-valued solutions of the
equation (\ref{Higgs_eq}). The $x$-independent solution of  (\ref{Higgs_eq})
solves the Duffing's-type equation
\[
 \ddot \psi +3 \dot \psi =\mu^2  \psi  -\lambda   \psi^ 3\,,
\]
which describes the motion of a mechanical system in a twin-well potential field.
Unlike  the equation in the Minkowski space-time, that is,  the equation
\begin{equation}
\label{1.1_Min}
\label{Higgs_eq_Min}
  \psi_{tt}  - \Delta  \psi =    \mu^2  \psi  -\lambda   \psi ^3\,,
\end{equation}
the equation (\ref{1.1}) has no other time-independent solution. For the  equation (\ref{1.1_Min}) the existence of a weak global solution in the energy space is known (see, e.g.,
\cite{G-V,G-V1989}).
The equation (\ref{1.1_Min})
for the Higgs scalar    field  in the Minkowski space-time has the time-independent flat solution
\begin{equation}
\label{stadysol}
\psi_M (x) = \frac{\mu }{\sqrt{\lambda} } \tanh \left(   \frac{\mu^2 }{2 } N\cdot (x-x_0) \right), \qquad N, x_0, x \in {\mathbb R}^3,
\end{equation}
where  $N $  is the unit vector.
The  solution (\ref{stadysol}), after Lorentz transformation, gives rise
to  a traveling solitary wave 
\[
\psi_M (x,t) =  \frac{\mu }{\sqrt{\lambda} } \tanh \left(   \frac{\mu^2 }{2 }  [ N \cdot (x-x_0) \pm v(t-t_0)]\frac{1}{\sqrt{1-v^2}}\right),
\qquad 
\]
 where  $ N, x_0, x \in {\mathbb R}^3$,
$t \geq t_0$, and $0<v<1$. The set of zeros of the solitary wave $\psi =\psi_M (x,t) $  is the moving boundary of the {\it wall}. 
\smallskip

A global  in time solvability of the Cauchy problem for equation (\ref{1.1})
 is not known, and the only  estimate for the lifespan is given by Theorem~0.1~\cite{ArXiv2017}.
The local solution exists
 for every smooth initial data. (See, e.g., \cite{Rendall_book}.)
The $C^2$ solution of the equation (\ref{Higgs_eq}) is unique and obeys the finite speed of the propagation
property. (See, e.g., \cite{Hormander_1997}.)
\smallskip

In order to make  our discussion   more transparent
we appeal to  the function \,$u = e^{\frac{3}{2}t}\psi$. For this  new unknown function \,$u=u(x,t) $,  the equation
(\ref{1.1})  takes the form of the semilinear
Klein-Gordon  equation  
\begin{equation}
\label{2.8b}
\label{2.8}
\label{2.2}
u_{tt} - e^{-2t} \bigtriangleup u  - M^2 u=  - \lambda e^{-3t} u^3,
\end{equation}
where a positive number \,$M $\, is defined as follows:
\[
M^2:=  \frac{9}{4} + \mu ^2>  0\,.
\]
The equation (\ref{2.2}) is the    equation with  imaginary mass. 
Next, we use the fundamental solution of the corresponding linear operator in order to reduce the Cauchy problem for the
semilinear equation to the integral equation and to define a weak solution.
We denote by $G$ the resolving operator of the problem
\[
u_{tt} - e^{-2t} \bigtriangleup u  - M^2 u=   f,   \quad u  (x,0) =  0 , \quad \partial_{t }u  (x,0 ) =0\,.
\]
Thus, $u=G[f]$. 
The operator $G$ is explicitly written   in \cite{Yag_Galst_CMP} for the case of the real mass.
The analytic continuation with respect to the parameter $M $  of this operator  allows us  also to use $G$
in the case of imaginary mass.  More precisely, for $M \geq 0$
 we define the operator \,$G$\, acting on \,$f(x,t) \in C^\infty ({\mb R}^3\times [0,\infty)) $ \,by (\ref{1.29small}), 
\begin{eqnarray*}
 G[f](x,t)
&  =  &
2   \int_{ 0}^{t} db
  \int_{ 0}^{ e^{-b}- e^{-t} } dr  \,  v(x,r ;b)    E(r,t;0,b;M)\,,
\end{eqnarray*}
where  the function
$v(x,t;b)$
is a solution to the Cauchy problem for the  wave equation (\ref{1.6c}).

Let $u_0=u_0(x,t)$ be a solution of the Cauchy problem
\begin{equation}
\label{3.5}
\partial_{t }^2 u_0 - e^{-2t} \bigtriangleup u_0   - M^2 u_0 =   0, \quad u_0 (x,0) = \varphi _0(x), \quad \partial_{t }u_0 (x,0 ) = \varphi _1(x )\,.
\end{equation}
Then any solution $u=u(x,t)$ of the equation (\ref{2.2}), which takes initial value $ u (x,0) = \varphi _0(x), \quad \partial_{t }u (x,0) = \varphi _1(x)$,
solves the integral equation
\begin{equation}
\label{2.7}
u(x,t) = u_0(x,t)- G[\lambda e^{-3\cdot } u^3](x,t)   \,.
\end{equation}
We use the last equation   to  define  a weak solution of the problem for the partial differential equation.
\begin{definition}
 If $u_0$ is a solution of the Cauchy problem (\ref{3.5}), then the solution  $u=u (x,t)$ of (\ref{2.7})
is said to be
{\it a weak solution} of the Cauchy problem for the equation (\ref{2.8}) with the initial conditions
$
u (x,0) = \varphi _0 (x)$,\,   $\partial_{t }u (x,0) = \varphi _1 (x) .
$
\end{definition}

It is suggested in \cite{CPDE2012} to measure a variation of the sign of the function $\psi   $ by the deviation from the H\"older inequality
\[
  \left| \int_{{\mathbb R}^n} u  (x) \, dx  \right|^{3}
\leq  C_{supp \, u}  \int_{{\mathbb R}^n}  |u   (x )|^3  \, dx
  \]
of the inequality between the
 integral of the function  and the self-interaction functional:
\[
  \left| \int_{{\mathbb R}^n} u  (x) \, dx  \right|^{3}
\leq \nu (u) \left| \int_{{\mathbb R}^n}  u ^3  (x )  \, dx \right|
  \]
provided that $\int_{{\mathbb R}^n}  u ^3  (x )  \, dx \not=0 $. For the
solutions with the initial data with supports in some bounded ball of radius $R$ due to finite speed of propagation the constant  $C_{supp \, u} $ depends of $R$ alone and is the same for all solutions. The constant $\nu (u)$ depends on function, but for the solution $u=u(x,t) $ of the equation (\ref{2.7}) it is regulated by the equation, that is, $\nu (u)=\nu (t)$ is a function of time universal for all functions. For the sign preserving global in time solutions
the rate of growth of the function $\nu (t)$ is restricted from below.

The   next definition  is a particular case of Definition~1.2~\cite{CPDE2012}.
Time $t$ is regarded as a parameter.
\begin{definition}
\label{D1.2}
The real valued-function $\psi  \in  C([0,\infty); L^1(  {\mb R}^3 )\cap L^3(  {\mb R}^3 ))$
is said to be asymptotically time-weighted   $L^3$-non-positive  (non-negative), if
there exist  number $C_\psi>0 $ and positive non-decreasing
function $\nu_\psi \in  C ([0,\infty)) $ such that with $\sigma =1$ ($\sigma =-1$)  one has
\[
  \left| \int_{{\mathbb R}^n} \psi  (x,t) \, dx  \right|^{3}
\leq -\sigma C_\psi \nu_\psi  (t) \int_{{\mathbb R}^n}   \psi ^3  (x,t)  \, dx
\quad \mbox{for all sufficiently large} \quad
 t .
  \]
\end{definition}
It is evident  that any sign preserving function $\psi  \in   L^3(  {\mb R}^3 ) $ with a compact support
satisfies the last inequality with $ \nu_\psi  (t) \equiv 1 $
and either $\sigma =1 $ or $\sigma =-1$, while $C_\psi^{ 1/2} $ is a measure of the support.

As a result of the finite speed of propagation property of the equation (\ref{1.1_Min}), any smooth global non-positive (non-negative)  solution $\psi = \psi (x,t) $ of   (\ref{1.1_Min}) with compactly supported initial data is
also asymptotically time-weighted   $L^3$-non-positive (non-negative) with the weight $\nu_\psi  (t)= (1+t)^6 $.

The following statement follows from Theorem~1.3~\cite{CPDE2012}.
Let $u=u(x,t) \in C([0,\infty);L^q(  {\mb R}^3 ))$, $2\leq q < \infty $, be a global solution
of the equation
\begin{eqnarray}
\label{210}
u(x,t) = u_0(x,t)- G\left[ \lambda    u^3 (y,\cdot ) \right](x,t)    \,.
\end{eqnarray}
where $ u_0(x,t)$ solves initial value problem  (\ref{3.5}) with $ \varphi _0, \varphi _1 \in C_0^\infty ({\mathbb R}^n)$ such that
\begin{eqnarray}
\label{21}
\sigma \left( M \int_{{\mathbb R}^n} \varphi _0(x )  \, dx + \int_{{\mathbb R}^n}  \varphi _1(x )  \, dx \right)>0\,.
\end{eqnarray}
 Assume also that the self-interaction functional satisfies
\[ 
\sigma  \int_{{\mathbb R}^n}    u^3(z, t)\, dz
\leq 0
\]
for all \, $t$\,   outside of the sufficiently small neighborhood of zero.
Then, the global solution $u=u(x,t)$ cannot be
an  asymptotically   time-weighted $L^3$-non-positive (-non-negative)  with the weight
$\nu_u =const>0$.

Thus, the the last statement shows that the continuous   
 global solution of the equation  (\ref{210}) cannot be negative  sign preserving  provided that
  it is generated by the function $u_0=u_0(x,t) $, which obeys (\ref{21}).
 Thus, {\it it takes positive value at some point, that is, it changes a sign.}

An application of the last theorem to the
 Higgs real-valued scalar field
equation (\ref{1.1})
with $\mu >0$ results in the following statement (see also Corollary~1.4~\cite{CPDE2012}). 
Let $\psi =\psi (x,t) \in C([0,\infty);L^q(  {\mb R}^3 ))$, $2\leq q < \infty $, be a global   weak solution
of the equation  (\ref{1.1}).
Assume also that
 the initial data of $\psi =\psi (x,t)$ satisfy
\begin{eqnarray}
\label{42}
\sigma \left( \left(\sqrt{9 + 4\mu ^2} +3\right)\int_{{\mathbb R}^3}\psi (x,0)\,dx+2\int_{{\mathbb R}^3}\partial _t\psi (x,0)\,dx \right) >0
\end{eqnarray}
with $\sigma =1$ ($\sigma =-1$), while
\begin{eqnarray*}
\sigma  \int_{{\mathbb R}^3}   \psi^3 (x,t)\, dx  \leq  0
\end{eqnarray*}
is fulfilled
 for all\, $t$\,  outside of the sufficiently small neighborhood of zero.

Then, the global solution $\psi =\psi (x,t)$ cannot be an  asymptotically   time-weighted $L^3$-non-positive (-non-negative)  solution  with the weight
$\nu_\psi   (t)=e^{a_\psi  t }t^{b_\psi } $, where
$ a_\psi <   \sqrt{9+4\mu ^2}- 3$, $b_\psi \in{\mathbb R} $.

For the solution $\psi =\psi (x,t) \in C^2( {\mb R}^3\times   [0,\infty))$ with the compactly supported smooth initial data $ \psi (x,0), \psi_t (x,0) \in C^\infty_0( {\mb R}^3 )$, the finite propagation speed property for
(\ref{1.1})
with $\mu >0$ implies that the solution has a  support in some cylinder $B_R \times [0,\infty) $, and consequently, if it is sign preserving, it is also
asymptotically   time-weighted $L^3$-non-positive (-non-negative)  solution  with the weight
$\nu_\psi   (t)\equiv 1 $. This contradicts  to the previous statement. Hence, {\it the global solution
with data satisfying (\ref{42}) and $\psi (x,0) \leq 0$
must take  positive value at some point and, consequently, must take zero value inside of some section $t=const>0$.}
It gives rise to the creation of a bubble.

\medskip

\section{Evolution of bubbles}
\label{S7}

Since an issue of the global solution for equation (\ref{1.1}) is not resolved, we present some simulation that shows evolution of the bubbles in time.

\begin{figure}
\caption{\label{fig:two-bubbles-init}Computational domain and first initial
data}
\centering{} 
\begin{tabular}{cc}
Computational domain  & Initial data along a diagonal line segment\tabularnewline
\includegraphics[height=4cm]{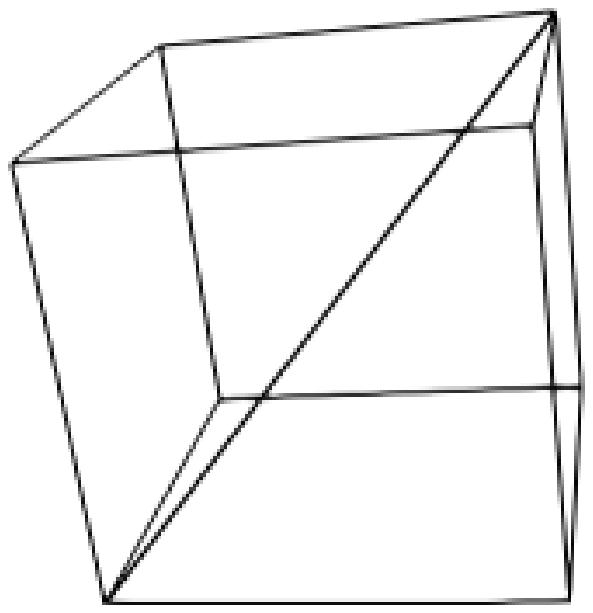} & \includegraphics[height=4cm]{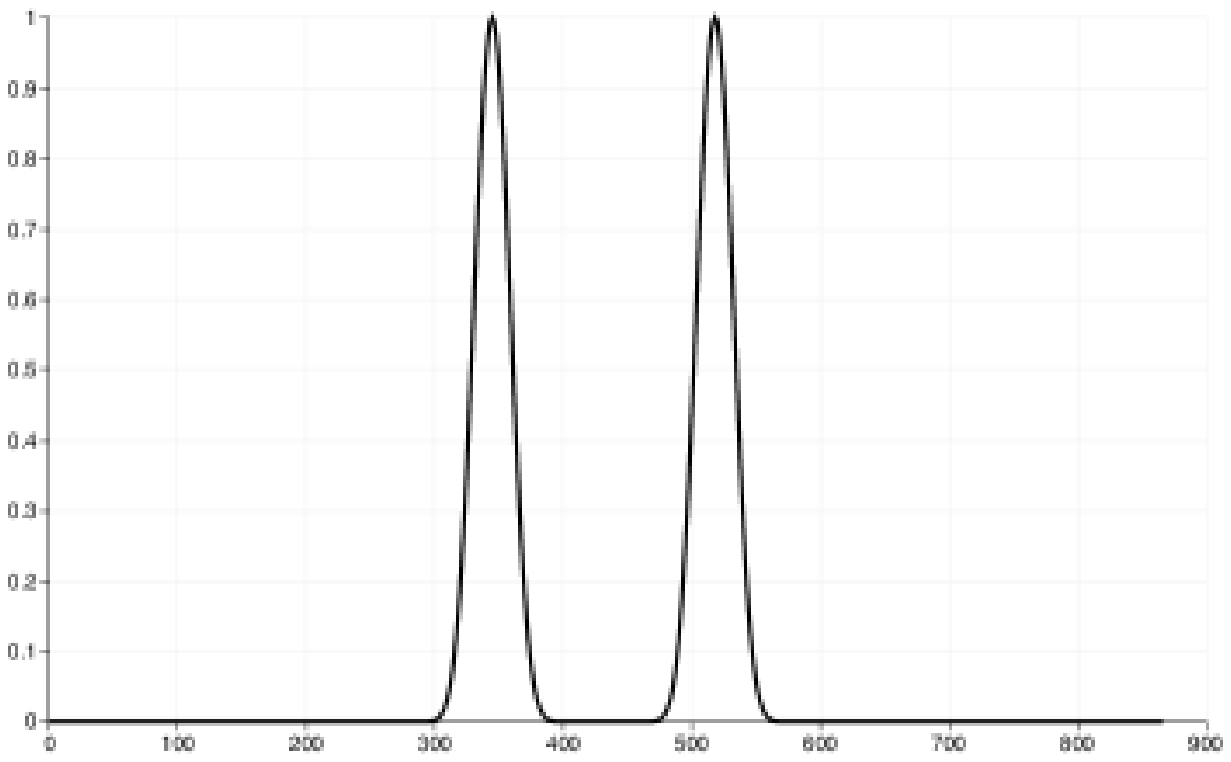}\tabularnewline
\end{tabular}
\end{figure}

\begin{figure}[h]
\caption{\label{fig:two-bubbles-2}Formation and interaction of two bubbles}
\centering{} 
\begin{tabular}{cc}
(a) 3D bubbles at $t=0.08$ & (b) 3D bubbles at $t=0.2$\tabularnewline
\includegraphics[height=3.5cm]{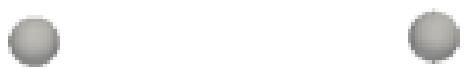} & \includegraphics[height=3.5cm]{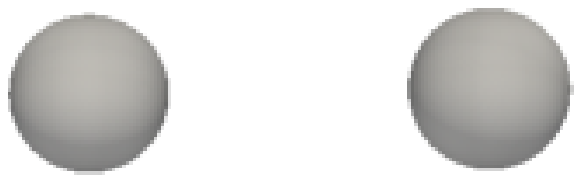}\tabularnewline
 & \tabularnewline
(c) 3D bubbles at $t=0.69$ & (d) 3D bubbles at $t=0.8$\tabularnewline
\includegraphics[height=3.5cm]{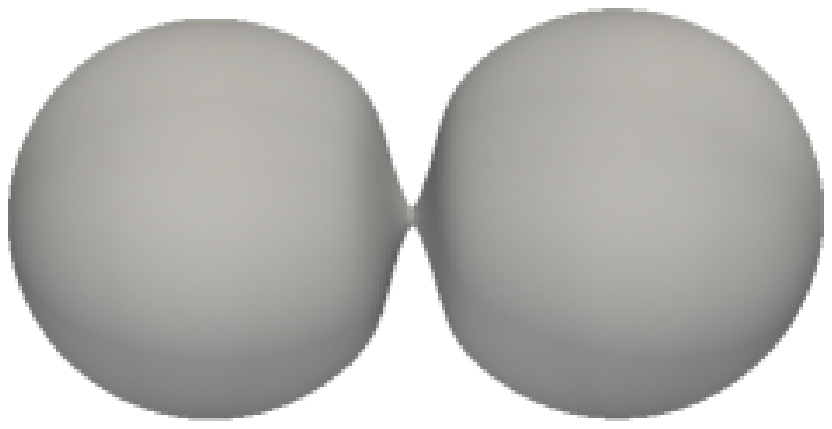} & \includegraphics[height=3.5cm]{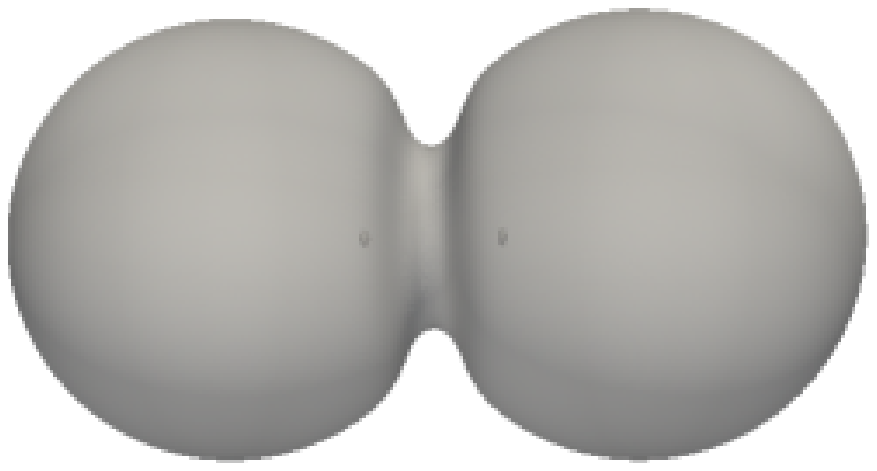}\tabularnewline
 & \tabularnewline
(e) 3D bubbles at $t=1$ & (f) 3D bubbles at $t=2$\tabularnewline
\includegraphics[height=3.5cm]{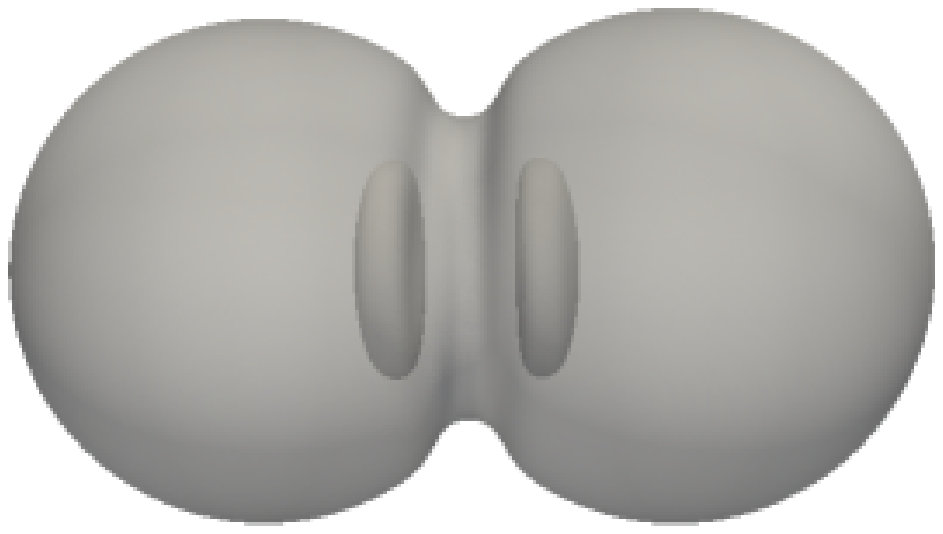} & \includegraphics[height=3.5cm]{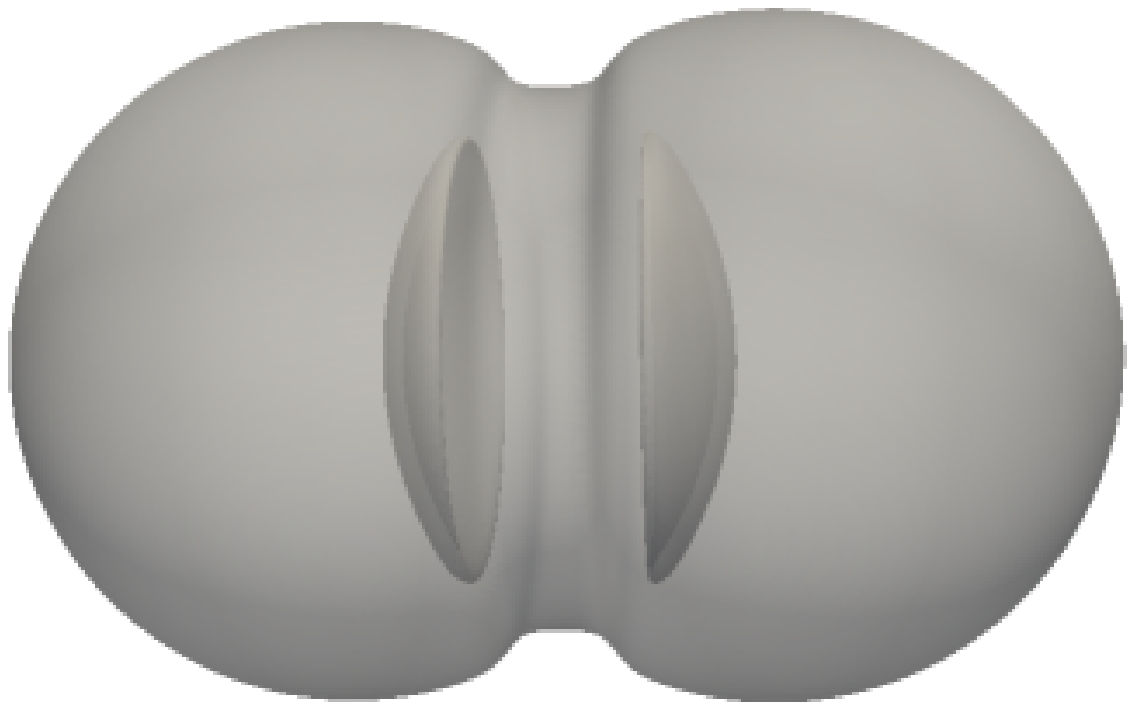}\tabularnewline
 & \tabularnewline
(g) 3D bubbles at $t=2.15$ & (h) 3D bubbles at $t=3$\tabularnewline
\includegraphics[height=3.5cm]{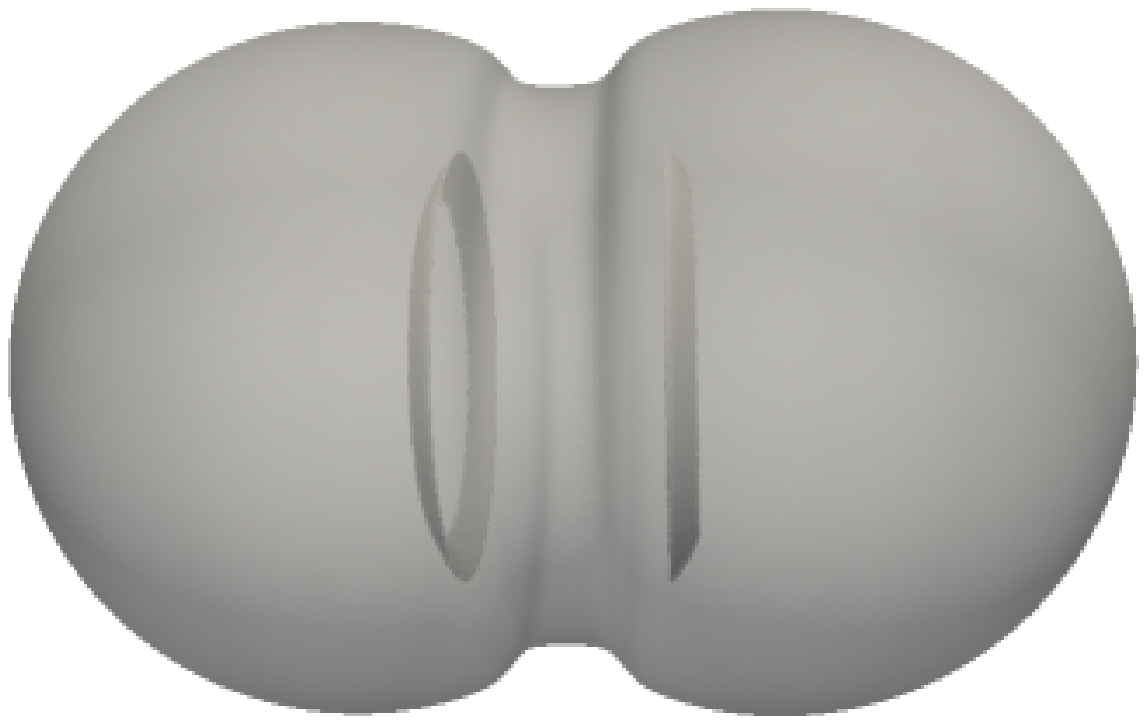} & \includegraphics[height=3.5cm]{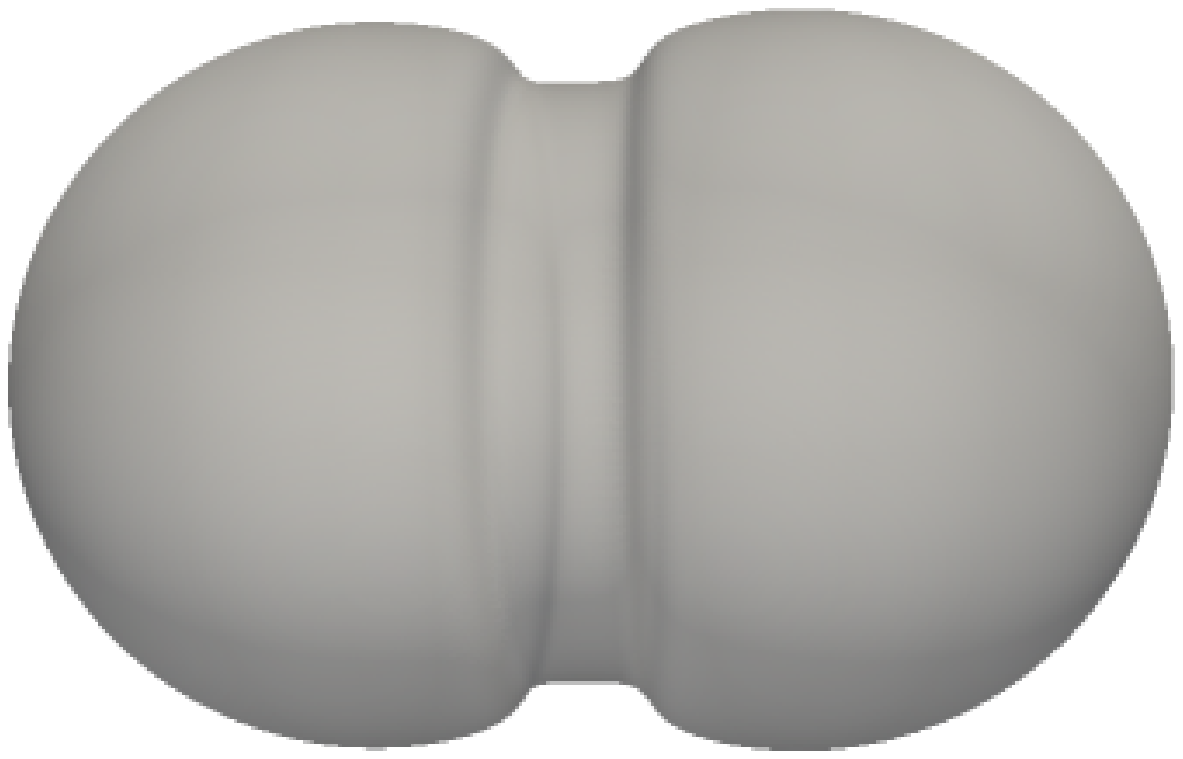}\tabularnewline
\end{tabular}
\end{figure}

Our numerical approach uses a fourth order finite difference method
in space \cite{finitediff} and an explicit fourth order Runge-Kutta
method in time \cite{RKbook} for the discretization of the Higgs
boson equation. The numerical code has been programmed using the Community
Edition of PGI CUDA Fortran \cite{PGI} on NVIDIA Tesla K40c GPU Accelerators.
The grid size in space was $n\times n\times n=501\times501\times501$,
resulting in a uniform spatial grid spacing of $\delta x_1=\delta x_2=\delta x_3=2\times10^{-2}$.
The time step $\delta t=10^{-3}$ ensured that the Courant\textendash Friedrichs\textendash Lewy
(CFL) condition \cite{Strang-Comp_Sci} for stability $\left({\displaystyle \left|\psi\right|<\frac{\delta x}{\sqrt{3}\delta t}\approx11.54}\right)$
was satisfied for all time. As first initial data $\psi_{0}$ we choose
the combination of two bell-shaped, infinitely smooth exponential
functions
\[
\psi_{0} (x )=B_{1} (x)+B_{2}(x) \qquad \forall x=(x_1,x_2,x_3)  \in\Omega,\label{eq:B1B2}, 
\]
where 
\begin{eqnarray*}
 B_{i}\left( {x}\right)=
\cases{ \exp\left(\frac{1}{R_{i}^{2}}-\frac{1}{R_{i}^{2}-\left| {x}-C_{i}\right|^{2}} \right) 
 \quad   if \quad \left|  {x}-C_{i}\right|<R_{i}, \cr 
0  \hspace{3.7cm}    if  \quad \left|  {x}-C_{i} \right|\geq R_{i}  }
\end{eqnarray*}
for $i=1,2$ with the center of the bell-shapes at $C_{1}=\left(0.4,0.4,0.4\right)$,
$C_{2}=\left(0.6,0.6,0.6\right)$, and the radii of the bell-shapes
$R_{1}=R_{2}=0.2$. Figure \ref{fig:two-bubbles-init} shows the computational
domain with a diagonal line segment and the line plot of the first
initial data $\psi_{0}\left(\vec{x}\right)$ along that line segment.
Note that the initial data is nonnegative with a compact support.
The finite cone of influence \cite{Hormander_1997} enables us to
use zero boundary conditions on the unit box $\Omega=\left(0,1\right)\times\left(0,1\right)\times\left(0,1\right)$
as computational domain, since the solution's domain of support stayed
inside the unit box. As second initial data $\psi_{1}$ we choose
a constant multiple of the first initial data
\[
\psi_{1} (x )   =-5\phi_{0}(x ) \qquad\forall   x=(x_1,x_2,x_3)  \in \Omega.\label{eq:5phi0}
\]
The parameter values are $\lambda=\mu^{2}=0.1$. Initially there is
no bubble present. Figure \ref{fig:two-bubbles-2} shows the formation
and interactions of bubbles. After the two bubbles form around time
$t=0.08$, their size grows continuously. Around time $t=0.69$ the
two bubbles touch, and from that time on they are attached to each
other. At time $t=0.8$ (shown on part (d) of Figure \ref{fig:two-bubbles-2})
an additional tiny bubble forms inside each of the now merged bubbles.
These additional bubbles grow (part (e) of Figure \ref{fig:two-bubbles-2}
at time $t=1$); then they flatten and become concave (part (f) of
Figure \ref{fig:two-bubbles-2} at time $t=2$). Later hole forms
in them and they become toroidal (part (g) of Figure \ref{fig:two-bubbles-2}
at time $t=2.15$), and finally they disappear (part (h) of Figure
\ref{fig:two-bubbles-2} at time $t=3$). The growth of the larger
outer bubble exponentially slows down and it does not seem to change
shape after time $t=3$.

\section*{Acknowledgments}

The authors acknowledge the Texas Advanced Computing Center 
at The University of Texas at Austin for providing high performance
computing and visualization resources that have contributed to the
research results reported within this paper. URL:  {http://www.tacc.utexas.edu}.
We also gratefully acknowledge the support of NVIDIA Corporation with
the donation of the Tesla K40 GPU used for this research.
K.Y. was supported by  University of Texas Rio Grande Valley College of Sciences 2016-17  Research Enhancement Seed Grant.


\begin{thebibliography}{00}




 


 
\bibitem{B-E}
H.~Bateman, ~A.~Erdelyi,~ {Higher  Transcendental  Functions}. {vol. 1,2}, New York:~McGraw-Hill, 1953.  

\bibitem{B-F-K-L}
A.~Bers,  R.~Fox,    C.~G~Kuper,    S.~G.~Lipson,  The impossibility of free tachyons, in Relativity and Gravitation, eds. C. G. Kuper and Asher Peres, New York, Gordon and Breach Science Publishers,    1971,     41--46.


\bibitem{Coleman}
S.~Coleman,~ {Aspects of Symmetry: Selected Erice Lectures}.  Cambridge University Press, 1985.  

\bibitem{RKbook}
K.~Dekker,  J.G.~Verwer,~ {Stability of Runge-Kutta methods for Stiff Nonlinear Differential Equations.}, North-Holland, Amsterdam:~Elsevier Science Ltd., 1984.




\bibitem{Epstein-Moschella}
H.~Epstein,  U.~Moschella,  de Sitter tachyons and related topics. Comm. Math. Phys.  {  336}~(1)  ~(2015)  381--430  

\bibitem{NARW2017}
A.~Galstian,  K.~Yagdjian,  Global in time existence of self-interacting scalar field in de Sitter
space-times. Nonlinear Analysis: Real World Applications {  34}~(2017)  110--139.

\bibitem{G-V}
 J.~Ginibre,, G.~Velo,~ The global Cauchy problem for the nonlinear Klein-Gordon equation.   {Math. Z.}  189, no. 4: (1985) 487--505.  

\bibitem{G-V1989}
J.~Ginibre,  G.~Velo,~ The global Cauchy problem for the nonlinear Klein-Gordon equation. II.   
\emph{Ann. Inst. H. Poincar{\'e} Anal. Non Lin{\'e}aire}.   6, no. 1 (1989) 15--35. 


\bibitem{Hawking}
S.~W.~Hawking,   G.~F.~R.~Ellis, {The large scale structure of space-time}.
 Cambridge Monographs on Mathematical Physics, No. 1. London-New York:~Cambridge University Press, 1973.  

\bibitem{Higgs}
P.W.~Higgs,  Broken symmetries and the masses of gauge bosons.~ {Phys. Rev. Lett.}   13, no.~16 (1964) 508--509. 

\bibitem{Hormander_1997}
L.~H\"ormander, ~ {Lectures on nonlinear hyperbolic differential equations}. 
Berlin:~ Springer-Verlag, 1997.  


\bibitem{finitediff}
H.B.~Keller,  V.~Pereyra,  Symbolic generation of finite difference formulas. Math. Comp.  { 32}~(144)  (1978) 955--971.  



\bibitem{Klainerman1980}
S.~Klainerman, Global existence for nonlinear wave equations. Comm. Pure Appl.
Math. 33 , no. 1  (1980) 43--101.


\bibitem{Lee-Wick}
T.D.~Lee,~ G.C.~Wick, ~Vacuum stability and Vacuum Excitation in Spin-0 Field. ~ {Phys. Rev. D}   9 (8)  (1974)   2291--2316. 

\bibitem{Linde}
A.~Linde,~ {Particle Physics and Inflationary Cosmology}. Harwood, Chur,
Switzerland, 1990. 



\bibitem{Moller}   
C.~M{\o}ller,~  {The theory of relativity}. Oxford: Clarendon Press, 1952.   

\bibitem{Protter}
M.H. Protter, H.F. Weinberger, ~{Maximum Principles in Differential Equations}.   by Springer-Verlag New York Inc., 1984

\bibitem{PGI}
NVIDIA Corporation: PGI CUDA Fortran Compiler. ~(2017) 







\bibitem{Rendall_book}
A.~Rendall, ~ {Partial differential equations in general relativity}. 
Oxford Graduate Texts in Mathematics, 16,  Oxford: Oxford University Press, 2008.  
 
\bibitem{Sather}
D.~Sather,  A maximum property of Cauchy's problem for the wave operator. Arch. Rational Mech. Anal.  21  (1966) 303--309.  
 
\bibitem{Shatah} 
J.~Shatah, M.~Struwe, {Geometric wave equations}. 
Courant Lect. Notes   Math., 2. New York Univ., 
New York: Courant Inst.   Math. Sci., 1998.     

 

\bibitem{Sonego} 
S. Sonego,    V.~Faraoni, Huygens' principle and characteristic propagation property for waves in curved space-times. J. Math. Phys.  33  ,  no. 2  (1992) 625--632.
 
 
\bibitem{Strang-Comp_Sci}
G.~Strang,  ~ {Computational science and engineering}. Wellesley, MA:~Wellesley-Cambridge Press, 2007.


\bibitem{Speck} 
J.~Speck,   
Finite-time degeneration of hyperbolicity without blowup for quasilinear wave equations, Analysis{\&}PDE, 
  10, no. 8  (2017) 2001--2030.  


 
 
 
\bibitem{Vasy_2010}  
A.~Vasy,~The wave equation on asymptotically de Sitter-like spaces. ~ {Adv. Math.}  223, no. 1  (2010) 49--97.
 
 
\bibitem{Voronov} 
N.A.~Voronov,   A.L.~Dyshko,~,  N.B.~Konyukhova,~On the Stability of a Self-Similar Spherical Bubble of a Scalar Higgs Field in de Sitter Space. ~ {Physics of Atomic Nuclei}    68, no. 7: (2005) 1218--1226.


\bibitem{Weinstein}
A.~Weinstein,  ~ On a Cauchy problem with subharmonic initial values. Ann. Mat. Pura Appl. (4)  43  (1957)  325--340.

\bibitem{Yag2005} 
K.~Yagdjian,~Global existence in the Cauchy problem for nonlinear wave equations with variable speed of propagation,  New trends in the theory of hyperbolic equations,   Oper. Theory Adv. Appl., 159, Birkhäuser, Basel, (2005)  301--385.


\bibitem{Yag_Galst_CMP} 
K.~Yagdjian, A.~Galstian,~ Fundamental Solutions for the Klein-Gordon Equation  in  de~Sitter space-time, 
~\emph{Comm. Math. Phys.} 285 (2009) 293--344.  

\bibitem{yagdjian_DCDS}
K.~Yagdjian, The semilinear Klein-Gordon equation in de Sitter space-time, Discrete Contin. Dyn. Syst. Ser. S {  2}~(3) (2009) 679--696.  


\bibitem{CPDE2012}
K.~Yagdjian, ~{On the global solutions of the Higgs boson equation,} {Comm. Partial Differential Equations} ~{\bf 37}~(3) ~(2012)  447--478. 

\bibitem{JMAA_2012}
K.~Yagdjian, ~Global existence of the scalar field in de Sitter space-time, J. Math. Anal. Appl.  {  396}~(1)   (2012) 323--344. 

\bibitem{JMP2013}
K.~Yagdjian, ~{Huygens' Principle  for the Klein-Gordon equation in the de~Sitter space-time,}  {J. Math. Phys.}  ~{ 54}, no. 9  (2013)  091503.

\bibitem{MN}
K.~Yagdjian, Integral transform approach to solving Klein-Gordon equation with variable coefficients,   Mathematische Nachrichten,
{\bf 288}~(17-18)  (2015)  2129-2152.

\bibitem{ArXiv2017}
K.~Yagdjian, Global existence of the  self-interacting scalar field in  the de~Sitter universe, arXiv:1706.07703 


\end{thebibliography}
\end{document}